\documentclass[10pt,a4paper]{article}
\usepackage{graphicx}
\usepackage{mathrsfs}
\usepackage{color}
\usepackage[latin1]{inputenc}

\usepackage{libertine}
\usepackage[T1]{fontenc}
\usepackage{lipsum}
\usepackage{titlesec}
\usepackage[colorlinks=true,linkcolor=magenta,citecolor=blue]{hyperref}

\usepackage[cmex10]{amsmath}
\usepackage{stmaryrd}
\usepackage{color}
\usepackage{cmll}
\usepackage{amsmath,amssymb}
\usepackage[sloped]{fourier}
\usepackage{amsthm}
\usepackage{esint}
\usepackage{psfrag}

\author{S\'ebastien Alvarez}
\date{}
\title{Harmonic measures and the foliated geodesic flow for foliations with negatively curved leaves}

\begin{document}

\newtheorem*{mtheorem}{\bf Main Theorem}\newtheorem{corollary}{Theorem}
\newtheorem*{mcoro}{\bf Corollary}
\newcounter{theorem}[section]
\newtheorem{exemple}{\bf Exemple \rm}
\newtheorem{exercice}{\bf Exercice \rm}
\newtheorem{conj}[theorem]{\bf Conjecture}
\newtheorem{defi}[theorem]{\bf Definition}
\newtheorem{lemma}[theorem]{\bf Lemma}
\newtheorem{proposition}[theorem]{\bf Proposition}
\newtheorem{coro}[theorem]{\bf Corollary}
\newtheorem{theorem}[theorem]{\bf Theorem}
\newtheorem{rem}[theorem]{\bf Remark}
\newtheorem{ques}[theorem]{\bf Question}
\newtheorem{propr}[theorem]{\bf Property}
\newtheorem{question}{\bf Question}
\def\bp{\noindent{\it Proof. }}
\def\ep{\noindent{\hfill $\fbox{\,}$}\medskip\newline}
\renewcommand{\theequation}{\arabic{section}.\arabic{equation}}
\renewcommand{\thetheorem}{\arabic{section}.\arabic{theorem}}
\newcommand{\eps}{\varepsilon}
\newcommand{\disp}[1]{\displaystyle{\mathstrut#1}}
\newcommand{\fra}[2]{\displaystyle\frac{\mathstrut#1}{\mathstrut#2}}
\newcommand{\dif}{{\rm Diff}}
\newcommand{\homeo}{{\rm Homeo}}
\newcommand{\Per}{{\rm Per}}
\newcommand{\Fix}{{\rm Fix}}
\newcommand{\A}{\mathcal A}
\newcommand{\Z}{\mathbb Z}
\newcommand{\Q}{\mathbb Q}
\newcommand{\R}{\mathbb R}
\newcommand{\C}{\mathbb C}
\newcommand{\N}{\mathbb N}
\newcommand{\T}{\mathbb T}
\newcommand{\U}{\mathbb U}
\newcommand{\D}{\mathbb D}
\newcommand{\PP}{\mathbb P}
\newcommand{\Sp}{\mathbb S}
\newcommand{\K}{\mathbb K}
\newcommand{\car}{\mathbf 1}
\newcommand{\g}{\mathfrak g}
\newcommand{\gs}{\mathfrak s}
\newcommand{\h}{\mathfrak h}
\newcommand{\rr}{\mathfrak r}
\newcommand{\s}{\sigma}
\newcommand{\fhi}{\varphi}
\newcommand{\ffhi}{\tilde{\varphi}}
\newcommand{\moins}{\setminus}
\newcommand{\ds}{\subset}
\newcommand{\W}{\mathcal W}
\newcommand{\WW}{\widetilde{W}}
\newcommand{\F}{\mathcal F}
\newcommand{\G}{\mathcal G}
\newcommand{\CC}{\mathcal C}
\newcommand{\RR}{\mathcal R}
\newcommand{\DD}{\mathcal D}
\newcommand{\M}{\mathcal M}
\newcommand{\B}{\mathcal B}
\newcommand{\cS}{\mathcal S}
\newcommand{\HH}{\mathcal H}
\newcommand{\Hyp}{\mathbb H}
\newcommand{\UU}{\mathcal U}
\newcommand{\Pp}{\mathcal P}
\newcommand{\QQ}{\mathcal Q}
\newcommand{\E}{\mathcal E}
\newcommand{\GG}{\Gamma}
\newcommand{\LL}{\mathcal L}
\newcommand{\KK}{\mathcal K}
\newcommand{\TT}{\mathcal T}
\newcommand{\X}{\mathcal X}
\newcommand{\Y}{\mathcal Y}
\newcommand{\ZZ}{\mathcal Z}
\newcommand{\bE}{\overline{E}}
\newcommand{\bF}{\overline{F}}
\newcommand{\wF}{\widetilde{F}}
\newcommand{\hcF}{\widehat{\mathcal F}}
\newcommand{\bW}{\overline{W}}
\newcommand{\bcW}{\overline{\mathcal W}}
\newcommand{\tL}{\widetilde{L}}
\newcommand{\diam}{{\rm diam}}
\newcommand{\diag}{{\rm diag}}
\newcommand{\dom}{{\rm dom}}
\newcommand{\Jac}{{\rm Jac}}
\newcommand{\Aff}{{\rm Aff}}
\newcommand{\Plong}{{\rm Plong}}
\newcommand{\Tr}{{\rm Tr}}
\newcommand{\Conv}{{\rm Conv}}
\newcommand{\Ext}{{\rm Ext}}
\newcommand{\Spec}{{\rm Sp}}
\newcommand{\Isom}{{\rm Isom}\,}
\newcommand{\Supp}{{\rm Supp}\,}
\newcommand{\Grass}{{\rm Grass}}
\newcommand{\Hold}{{\rm H\ddot{o}ld}}
\newcommand{\Ad}{{\rm Ad}}
\newcommand{\ad}{{\rm ad}}
\newcommand{\Aut}{{\rm Aut}}
\newcommand{\End}{{\rm End}}
\newcommand{\Leb}{{\rm Leb}}
\newcommand{\Int}{{\rm Int}}
\newcommand{\cc}{{\rm cc}}
\newcommand{\grad}{{\rm grad}}
\newcommand{\proj}{{\rm proj}}
\newcommand{\mass}{{\rm mass}}
\newcommand{\dive}{{\rm div}}
\newcommand{\dist}{{\rm dist}}
\newcommand{\im}{{\rm Im}}
\newcommand{\re}{{\rm Re}}
\newcommand{\codim}{{\rm codim}}
\newcommand{\Map}{\longmapsto}
\newcommand{\vide}{\emptyset}
\newcommand{\tr}{\pitchfork}
\newcommand{\ssl}{\mathfrak{sl}}

\newenvironment{demo}{\noindent{\textbf{Proof.}}}{\quad \hfill $\square$}
\newenvironment{pdemo}{\noindent{\textbf{Proof of the proposition.}}}{\quad \hfill $\square$}
\newenvironment{IDdemo}{\noindent{\textbf{Idea of proof.}}}{\quad \hfill $\square$}

\def\to{\mathop{\rightarrow}}
\def\act{\mathop{\curvearrowright}}
\def\To{\mathop{\longrightarrow}}
\def\Sup{\mathop{\rm Sup}}
\def\Max{\mathop{\rm Max}}
\def\Inf{\mathop{\rm Inf}}
\def\Min{\mathop{\rm Min}}
\def\lims{\mathop{\overline{\rm lim}}}
\def\limi{\mathop{\underline{\rm lim}}}
\def\egal{\mathop{=}}
\def\dans{\mathop{\subset}}
\def\surj{\mathop{\twoheadrightarrow}}

\def\h#1#2#3{ {\rm{hol}}^{#1}_{#2\rightarrow#3}}

\maketitle

\begin{abstract}
In this paper we define a notion of Gibbs measure for the geodesic flow tangent to a foliation with negatively curved leaves and associated to a particular potential $H$. We prove that there is a canonical bijective correspondence between these measures and Garnett's harmonic measures.
\end{abstract}

\section{Introduction}

This paper is the first of a series of three papers in which we study a notion of Gibbs measures for the geodesic flow tangent to the leaves of a foliation with negatively curved leaves \cite{Al1,Al2}.

In general, the holonomy pseudogroup of a foliation doesn't preserve a transverse measure. In order to overcome this difficulty Garnett developed in \cite{Gar} a notion of \emph{harmonic measures}. These are the measures which are invariant by heat diffusion along the leaves.

A link with the \emph{foliated geodesic flow} emerged in the work of Mart\'inez, Bakhtin-Mart\'inez. Consider $(M,\F)$ a compact foliated space whose leaves are \emph{hyperbolic Riemann surfaces}. It follows from \cite{BMar,M} that there is a canonical bijective correspondence between the harmonic measures for $\F$ and the measures on the unit tangent bundle $T^1\F$ which are invariant by the joint action of the foliated geodesic and unstable horocyclic flows.

The goal of this paper is to explain what happens when the leaves are only negatively curved manifolds. In that case, there is no horocyclic flow anymore and even if we were able to lift harmonic measures to the unit tangent bundle there is no reason why they should be invariant by the foliated geodesic flow.

The idea, aldready present in \cite{L2}, is that in negative curvature the Brownian paths are escorted by geodesics. In particular, almost every Brownian path converges to some point at infinity. It is therefore interesting to define a measure which is invariant by the foliated geodesic flow and which describes precisely the ergodic behaviour of geodesic escorts of the Brownian paths tangent to the leaves. The aim of the present work is to show that such measures exist and have the Gibbs property.

\paragraph{Harmonic and Gibbs measures.} Let $N$ be a complete, connected and simply connected manifold whose sectional curvature is pinched between two negative constants. The \emph{harmonic class} on the sphere at infinity $N(\infty)$ is the class of the hitting distributions at infinity of Brownian paths. This measure class can be projected down onto unstable horospheres along the geodesics, thus defining the harmonic class on horospheres. This harmonic class on the unstable horospheres is by definition preserved by the geodesic flow. Thus, when $(M,\F)$ is a compact foliated space with negatively curved leaves one can consider probability measures on $T^1\F$ which are invariant by the foliated geodesic flow, denoted by $G_t$, and whose conditional measures in the unstable manifolds lie in the harmonic class. Such a measure is called an $H$-\emph{Gibbs measure}. It is associated to the following potential:

$$H(v)=\left.\frac{d}{dt}\right|_{t=0}\log\,k_{L(v)}(c_v(0),c_v(t);c_v(-\infty)),$$
where $L(v)$ is the leaf $v$ is tangent to, $k_{L(v)}$ is the Poisson kernel on $\widetilde{L(v)}$ and $c_v$ is the geodesic directed by $v$. The main result of this paper is then:

\begin{mtheorem}
Let $(M,\F)$ be a compact foliated space whose leaves are of class $C^{\infty}$. Let $\TT$ be a complete transversal to $\F$. We endow each leaf with a $C^{\infty}$ Riemannian metric which varies continuously in the $C^{\infty}$-topology. Assume that the leaves have negative sectional curvature. Then:
\begin{itemize}
\item for any harmonic measure for $\F$, there is a unique $H$-Gibbs measure on $T^1\F$ for the foliated geodesic flow $G_t$ inducing the same measure on $\TT$.

\item reciprocally, for any $H$-Gibbs measure for $G_t$, there is a unique harmonic measure for $\F$ inducing the same measure on $\TT$.
\end{itemize}
\end{mtheorem}

When all the leaves are of constant negative curvature, the harmonic class in horospheres coincides with the Lebesgue class and the notion of $H$-Gibbs measures coincide with the usual notion of Gibbs $u$-states of \cite{BDV} (i.e. invariant measures with Lebesgue disintegration in the unstable manifolds). Moreover, when the leaves are $2$-dimensional these Gibbs $u$-states are exactly the measures invariant by both the foliated geodesic and unstable horocyclic flows. The following corollary is thus a generalization of the work of Mart\'inez, Bakhtin-Mart\'inez because it treats the case of foliations by higher dimensional manifolds. Note that in \cite{CM}, this corollary has even been generalized in the case of foliations by locally symmetric spaces.

\begin{mcoro}
Let $(M,\F)$ be a compact foliated space whose leaves are hyperbolic manifolds. Then a measure is harmonic if and only if it is the projection of a Gibbs $u$-state for the foliated geodesic flow $G_t$. Moreover, two different Gibbs $u$-states for $G_t$ project down to different harmonic measures for $\F$.
\end{mcoro}

\paragraph{The unrolling argument.} The main purpose of this work is to give a short argument in order to prove the main Theorem. This argument is based on the integral representation of harmonic functions in negative curvature. We summarize it below:
\begin{itemize}
\item first, we will see in Proposition \ref{liftharmonicmeasures} (see also the proposition below for a simple illustration) how to lift canonically to $T^1\F$ a harmonic measure for a lamination $\F$ with negatively curved leaves by ``unrolling'' its local harmonic densities;
\item second, we will see in Proposition \ref{reparamharmhgibbs} how to reparametrize this measure in the unstable manifolds so as to obtain an invariant measure by the foliated geodesic flow which is an $H$-Gibbs measure.
\end{itemize}

Let us briefly explain what we mean by ``unrolling'' the densities and treat a very simple example. Consider the upper half plane $\Hyp$ endowed with the Poincaré metric $ds^2=y^{-2}(dx^2+dy^2)$. One has a trivialization $T^1\Hyp=\Hyp\times S_{\infty}$ (here, $S_{\infty}=\R\PP^1$ denotes the circle at infinity) obtained by identifying $\Hyp\times\{\xi\}$ with the corresponding center unstable leaf of the geodesic flow. The following proposition was suggested by Marco Brunella.

\begin{proposition}
A Borel measure on $\Hyp$ is harmonic if and only if it is the projection of a Borel measure on $T^1\Hyp$ which is invariant by the joint action of the geodesic and unstable horocyclic flows. Moreover, two different measures on $T^1\Hyp$ invariant by the joint action of the geodesic and horocyclic flows project down onto different harmonic measures on $\Hyp$.
\end{proposition}

\begin{IDdemo}
On the one hand, a harmonic measure of $\Hyp$ writes as $m=h\Leb$ (here $\Leb$ denotes the hyperbolic area) where $h:\Hyp\to(0,\infty)$ is a harmonic function: $\Delta h=0$. But such a function possesses an integral representation: $h(z)=\int_{S_{\infty}}k(z;\xi)\,d\eta(\xi)$, where $k$ is the usual \emph{Poisson kernel} and $\eta$ is a (uniquely determined) finite Borel measure on $S_{\infty}$.

On the other hand, we associate to each center unstable leaf $\Hyp\times\{\xi\}$ the measure $\mu_{\xi}=k(z,\xi)\Leb(z)$. Up to a multiplicative constant this is the unique measure on this leaf which is invariant by the joint action of the geodesic and unstable horocyclic flows.

In order to see this, identify $\Hyp\times\{\infty\}$ with the group $\Aff<PSL_2(\R)$ constituted of lower triangular matrices. The center unstable leaves then correspond to the orbits of $\Aff$ on $PSL_2(\R)$ by right translations. In particular, the joint action of the flows on $\Hyp\times\{\infty\}$ corresponds to the action of $\Aff$ on itself by right translations. The unique right Haar measure reads in coordinates as $dxdy/y=y\Leb$. But we have $k(x+\mathbf{i}y;\infty)=y$. Thus by composing by the homography $z\mapsto-1/(z-\xi)$ we see that on $\Hyp\times\{\xi\}$ the unique measure invariant by the joint action of the flows is given by $\mu_{\xi}$.

Now integrate the measures $\mu_{\xi}$ against $\eta$. We obtain a measure $\mu$ on $T^1\Hyp$ which projects down to $m$ and is invariant by both the geodesic and unstable horocyclic flows. It can be seen that any measure invariant by this joint action may be written this way.

Two different measures invariant by the two flows project down onto different measures because the density of the projection is determined by the measure $\eta$ on $S_{\infty}$: the integral representation of positive harmonic functions is unique.

We obtained $\mu$ by unrolling the density $h$: we call it the \emph{canonical lift} of $m$. This operation gives a canonical bijective correspondence between harmonic measures on $\Hyp$ and measures on $T^1\Hyp$ which are invariant by the joint action of the geodesic and unstable horocyclic flows.
\end{IDdemo}\\

In the case of laminations, the idea is to repeat this argument leaf by leaf, while being careful with the cocycle relations imposed by the holonomy.

The unrolling argument will be used in the context of Gibbs u-states for the foliated geodesic flow in \cite{Al1} and of Gibbs measures for the foliated bundles associated to any potential in \cite{Al2}.

\section{Harmonic measures and the foliated geodesic flow}

\subsection{Foliated spaces and holonomy}

\paragraph{Laminations.} Let $M$ be a compact metrizable space. A \emph{lamination} $\F$ of $M$ is the data of a finite foliated atlas $\A=(U_i,\phi_i)_{i\in I}$. Such an atlas consists of an open cover $(U_i)_{i\in I}$ of $M$ and of homeomorphisms $\phi_i:U_i\to P_i\times T_i$ where $P_i$ is a cube of dimension $d$ (this integer is called the \emph{dimension} of the lamination) and $T_i$ is a topological space, such that the changes of charts are of the form:
\begin{equation}
\label{Eq:changecharts}
\phi_j\circ\phi_i^{-1}(z,x)=(\zeta_{ij}(z,x),\tau_{ij}(x)),
\end{equation}
where $(z,x)\in\phi_i(U_i\cap U_j)$, $\zeta_{ij}(.,x)$ being a $C^{\infty}$-diffeomorphism of $P_i$ which varies continuously in the $C^{\infty}$-topology with respect to $x$. We shall say that $(M,\F)$ is a \emph{foliated space}.

The sets $\phi_i^{-1}(P_i\times\{x\})$ are called \emph{plaques}. Since the changes of charts have this very particular form, we can glue the plaques together so as to obtain a partition of $M$ by manifolds called the \emph{leaves} of $\F$. We denote by $(P_i(x))_{x\in T_i,i\in I}$ the family of plaques of $\F$ and, abusively, we will always identify the transversals $T_i$ with their preimages $\phi_i^{-1}(\{z_i\}\times T_i)$ for some choice of $z_i$. Such a collection of spaces whose union meets each leaf is called a \emph{complete system of transversals}.

\paragraph{Holonomy.} Let $(M,\F)$ be a compact foliated space. Fix a foliated atlas $\A=(U_i,\phi_i)_{i\in I}$ of $M$ and assume that it is a \emph{good atlas} in the sense that:
\begin{enumerate}
\item if two plaques intersect each other, their intersection is connected;
\item if two charts $U_i$ and $U_j$ intersect each other, then $\overline{U_i\cap U_j}$ is included in a foliated chart: a plaque of $U_i$ intersects at most one plaque of $U_j$.
\end{enumerate}

We can define a \emph{complete transversal} $\TT$ as the union of all the $T_i$. The maps $\tau_{ij}$ generate a \emph{pseudogroup} $\Pp$ of local homeomorphisms of $\TT$, called the \emph{holonomy pseudogroup} of $\F$. Recall that we chose to identify $T_i$ and a local transversal $\phi_i^{-1}(\{z_i\}\times T_i)$ and to consider the points of $\TT$ as points of $M$. Note that this pseudogroup is symmetric since we have $\tau_{ij}^{-1}=\tau_{ji}$.

We will often talk about the \emph{holonomy map along a path}. If $c:[0,1]\to M$ is a path tangent to a leaf and $T_i$ and $T_j$ are local transversals to $\F$ containing $c(0)$ and $c(1)$, there exist two neighbourhoods $S_i\dans T_i$ and $S_j\dans T_j$ of $c(0)$ and $c(1)$ respectively such that we can send every $x\in S_i$ onto $\tau_c(x)\in S_j$ by sliding along the leaves of $\F$. More precisely, if we consider any chain of charts that cover $c$, say $U_{i_0},...,U_{i_n}$, then $\tau_c$ is defined as the composition $\tau_{i_{n-1}i_n}\circ...\circ\tau_{i_1 i_0}$. The germ at $c(0)$ of $\tau_c:S_i\to S_j$ does not depend on the choice of $c$ nor does it on the choice of the chain of charts. It only depends on the homotopy class of $c$.

Let us emphasize the fact that even if acts on the noncompact space $\TT$, the holonomy pseudogroup $\Pp$ is \emph{compactly generated} in the sense of Haefliger \cite{H}. That means that:
\begin{itemize}
\item the complete transversal $\TT$ is locally compact and contains a relatively compact open set $\TT_0$ which meets all orbits of $\Pp$;
\item there is a symmetric finite set $\cS=(\sigma_j)_{j\in J}$ of elements of $\Pp$ which generate the restriction $\Pp_{|\TT_0}$, and such that each $\sigma_j:S_j\to S_j'$ is the restriction of an element $\overline{\sigma}_j$ of $\Pp$ whose domain contains the adherence of $S_j$.
\end{itemize}

\subsection{Disintegration in the leaves of a foliation}

Since in general a lamination $\F$ of a compact space $M$ does not define a \emph{measurable partition} (see the terminology of Rokhlin \cite{Ro}), we can't a priori disintegrate a probability measure in the leaves of $\F$. However, because it is a locally trivial partition we  can do it locally and, given a probability measure on $M$, we can compare the class of conditional measures to a prescribed one, as explained below.

\begin{defi}
\label{disintegrationleaves}
Let $(\F,\lambda^{\F}_x)$ be a lamination of a compact space $M$, together with a measurable family of Borel measures on the leaves, such that for $x,y$ in the same leaf, $\lambda^{\F}_x=\lambda^{\F}_y$.

We say that $\mu$ has an absolutely continuous (resp. singular) disintegration with respect to $(\F,\lambda^{\F}_x)$ if its conditional measure in almost every plaque 
$P(x)$ is absolutely continuous (resp. singular) with respect to $\lambda^{\F}_x$.
\end{defi}

The most usual context is when $\lambda^{\F}_x$ is the Lebesgue measure on the leaves of $\F$. In that case, we say that $\mu$ has \emph{Lebesgue disintegration in the leaves of} $\F$.

\subsection{Leafwise metrics and foliated geodesic flow}

\paragraph{Leafwise metrics.} Let $(M,\F)$ be a compact foliated space. Assume that:
\begin{itemize}
\item each leaf $L$ possesses a $C^{\infty}$ Riemannian metric denoted by $g_L$;
\item the metric $g_L$ varies continuously with $L$ in the $C^{\infty}$-topology.
\end{itemize}
The data of all the $g_L$ gives what we call a leafwise metric for $\F$. As $M$ is compact, a leafwise metric gives \emph{uniformly bounded geometry} to the leaves of $\F$: their injectivity radii are uniformly bounded from below and their sectional curvatures are uniformly pinched.

\paragraph{The unit tangent bundle.} Assume that $(M,\F)$ is a compact foliated space with a leafwise metric. We can form the \emph{unit tangent bundle} of the foliation, which is denoted by $T^1\F$ and consists of the unitary vectors tangent to the leaves of $\F$. This is a compact space foliated by the $T^1L$, $L$ leaf of $\F$. This lamination will be denoted by $\widehat{\F}$. We assume from now on that the leaves are of dimension $\geq 2$ for in the case of a $1$-dimensional foliation, the space $T^1\F$ consists of two copies of $M$. It is important to notice that $\F$ and $\widehat{\F}$ have the same holonomy:

\begin{lemma}
\label{sameholonomy}
The holonomy pseudogroups of the two laminations $\F$ and $\widehat{\F}$ are the same.
\end{lemma}

\begin{demo}
Starting with a foliated atlas of $\F$, we can construct a foliated atlas of $\widehat{\F}$ by pulling it back by the fibration. Indeed, because of the uniformly bounded geometry of the leaves it is possible to choose the charts in $M$ in such a way that they trivialize the unit tangent bundle of $\F$. Hence, their preimages by the fibration form an atlas of $\widehat{\F}$ whose charts are trivially subfoliated by spheres of dimension $\dim\F-1$: the holonomy of this subfoliation is trivial. Thus, we did not add any holonomy by pulling back by the fibration.
\end{demo}

\paragraph{Remark.} In particular, a lamination $\F$ has a transverse measure invariant by holonomy if and only if the laminition $\hcF$ of $T^1\F$ has one.\\

We can endow $T^1\F$ with a leafwise metric by considering the foliated \emph{Sasaki metric}. The induced Riemannian distance is defined as follows. Locally, if $v$ and $w$ belong to the same leaf $T^1L$ the square of their distance is the sum of the square of the distance between the base points and of the square of their angle (computed with the parallel transport). This metric varies continuously in the $C^{\infty}$-topology with the transverse parameter because the Levi-Civita connexion only depends on the 1-jet of the metric (its coefficients are given by the Christoffel's symbols). In particular, it varies continuously with the transverse parameter in the $C^{\infty}$-topology.

\paragraph{The foliated geodesic flow.} A vector $v\in T^1\F$ directs a unique geodesic inside the leaf of $v$. We define $G_t(v)$ by flowing $v$ along this geodesic during a time $t$. This allows us to define the \emph{foliated geodesic flow} $G_t$. This flow is $C^{\infty}$ inside the leaves of $T^1\F$ and varies continuously in the $C^{\infty}$-topology with the transverse parameter. Indeed, the geodesics are solutions of the \emph{geodesic equation} which is a second order ODE whose coefficients are locally given by the Christoffel's symbols. Therefore, this equation as well as its solution is continuous with the metric in the $C^{\infty}$-topology.

\subsection{Harmonic measures and transverse measures}
In what follows, $(M,\F)$ is a compact foliated space with a leafwise metric.
\paragraph{Garnett's harmonic measures.} Given that each leaf is endowed with a smooth Riemannian structure, there is a Laplace-Beltrami operator on each leaf that we will denote by $\Delta_L$.

\begin{defi}
\label{laplace}
A probability measure $m$ on $M$ is said to be harmonic if it has Lebesgue disintegration in the leaves of $\F$ and if the local densities in the plaques are harmonic.
\end{defi}

Garnett showed in \cite{Gar} the existence of such measures. They describe the Brownian paths tangent to the leaves in the sense that they are invariant by leafwise heat diffusion.

The set of harmonic measures form a convex compact subset of probability measures in $M$. Extremal points are called \emph{ergodic} harmonic measures. More precisely, a harmonic measure is ergodic if it gives null or full measure to any saturated Borel set $X$. As in the case of invariant measures for flows or diffeomorphisms, an \emph{ergodic decomposition theorem} exists. Each harmonic measure may be written in a unique way as a barycenter of ergodic ones. See \cite{Gar} for the precise statements.

\paragraph{Cocycle relations.} To any harmonic measure can be associated a \emph{Radon-Nikodym cocycle} on a complete transversal in the following way. Let $\A=(U_i,\phi_i)_{i\in I}$ be a good foliated atlas. Write each chart as a union of plaques: $U_i=\bigcup_{x\in T_i} P_i(x)$, where $(T_i)_{i\in I}$ is a complete system of transversals. In a chart $U_i$, by the definition of harmonic measures, there is a measure $\nu_i$ on $T_i$ as well as a measurable family indexed by $x\in T_i$ of harmonic functions $h_i(.,x):P_i(x)\to (0,\infty)$ such that: $m_{|U_i}=h_i(z,x)\Leb_{|P_i(x)}(z)\, \nu_i(x).$

The evaluation of $m$ on the intersection of two charts $U_i$ and $U_j$ yields the following formula for $x\in T_i$ such that the plaques $P_i(x)$ and $P_j(\tau_{ij}(x))$ intersect and $z$ lying in the intersection:
\begin{equation}
\label{Eq:cocyclerelations1}
\frac{d[\tau_{ij}^{-1}\ast\nu_j]}{d\nu_i}(x)=\frac{h_i(z,x)}{h_j(z,\tau_{ij}(x))}.
\end{equation}

\paragraph{Transverse measure.} Consider the complete transversal defined as $\TT=\bigsqcup_{i\in I} T_i$. The holonomy pseudogroup $\Pp$ acts on $\TT$. Then, it is possible to consider the measure $\widehat{m}$ on $\TT$ which is equal to $\nu_i$ when restricted to $T_i$. We call this measure the \emph{measure induced} by $m$ on $\TT$.

Relation (\ref{Eq:cocyclerelations1}) shows that the measure $\widehat{m}$ is \emph{quasi-invariant} by the action of the holonomy pseudogroup: if $A\dans\TT$ lies in the domain of a holonomy map $\tau$ and $\widehat{m}(A)>0$, then $\widehat{m}(\tau(A))>0$.

\paragraph{Harmonic extension.} A priori the presence of holonomy in the leaves is an obstruction to extend harmonically the local densities of $m$. The following lemma due to Ghys (see \cite{Gh}) shows that for almost every leaf, this obstruction is in fact void.

\begin{lemma}
\label{noobstruction}
There exists a Borel set $\X_0\dans\TT$ which is full for $\widehat{m}$ and saturated by the action of the pseudogroup $\Pp$, such that for any $x\in\X_0$ and any $\gamma\in\Pp$ which fixes $x$, we have:
$$\frac{d[\gamma\ast\widehat{m}]}{d\widehat{m}}(x)=1.$$
\end{lemma}

In other terms, Ghys' Lemma shows that for all $i\in I$ and $m$-almost $z\in U_i$, if $z\in P_i(x)$ and $c_1$, $c_2$ are two path tangent to $L_z$ which both start at $z$ and end at a same point $z'\in U_j$, then:
$$\frac{d[\tau_{c_1}^{-1}\ast\nu_j]}{d\nu_i}(x)=\frac{d[\tau_{c_2}^{-1}\ast\nu_j]}{d\nu_i}(x).$$

Now it is  easy to show how to extend a local density of $m$ to a whole leaf:

\begin{lemma}
\label{harmonicextension}
Let $m$ be a harmonic measure for $\F$. Then there is a full and saturated Borel set $\X\dans M$ such that for any $y\in\X$, if $y\in P_{i_0}(x)\dans U_{i_0}$, for some $i_0\in I$ and $x\in T_{i_0}$, the following formula defines a positive harmonic function of $z\in L_y$:
\begin{equation}
\label{Eq:harmonicextension1}
H_x(z)=\frac{d[\tau_c^{-1}\ast\nu_i]}{d\nu_{i_0}}(x) h_i(z,\tau_c(x)),
\end{equation}
where $z\in U_i$ and $c$ is any path joining $y$ and $z$.
\end{lemma}

\section{Harmonic measures for foliations with negatively curved leaves}

\subsection{Brownian motion in negative curvature} Hereafter, $N$ denotes a complete, connected and simply connected Riemannian manifold of class $C^{\infty}$ whose sectional curvature is everywhere pinched between two negative constants $-b^2\leq-a^2$. The \emph{sphere at infinity} of $N$ is denoted by $N(\infty)$.

It is classical that the Laplace operator $\Delta$ generates a one-parameter semi group called the \emph{heat diffusion}, or \emph{Brownian diffusion} operator. This allows us for all $z$ to define naturally the \emph{Wiener measure} $W_z$ on the space $\Omega_z$ of continuous paths $\omega:[0,\infty)\to N$ with $\omega(0)=z$ (see \cite{Su} for more details). A \emph{Brownian path} starting at $z$ is a typical path for the $W_z$.

\paragraph{Harmonic class.}

Ledrappier showed in \cite{L2} Proposition 5 that a Brownian path on $N$ has a \emph{geodesic escort}: it stays at a sublinear distance to a geodesic ray. Note that he stated this proposition for Riemannian covers of closed and negatively curved manifolds but that the compactness hypothesis is not used in his proof: it only uses the pinching of the curvature. He attributes this argument to \cite{K,L1}. Note finally that his proof uses Prat's work about the evolution of polar coordinates along a typical Brownian path in pinched negative curvature \cite{Pr}: he proved that the growth of the radial coordinate of the time $t$ of a Brownian path is linear and that the angular coordinate has a limit as $t$ tends to infinity.

Sullivan considered in \cite{Su} the hitting distribution of Brownian paths at infinity. The family of measures on $N(\infty)$ he obtained, denoted by $(w_z)_{z\in N}$, satisfies the following properties:

\begin{itemize}
\item they are all in the same measure class, that we call \emph{harmonic class};
\item the harmonic class has no atoms and gives positive measure to any nonempty open set of $N(\infty)$. Better: when $z$ tends to $\xi$ in the cone topology $w_z$ tends to $\delta_{\xi}$.
\end{itemize}

\paragraph{Poisson kernel and integral representation of harmonic functions.} We define the \emph{Poisson kernel} of $N$ as the Radon-Nikodym cocycle:
\begin{equation}
\label{Eq:Poisson}
k(z_1,z_2;\xi)=\frac{dw_{z_2}}{dw_{z_1}}(\xi).
\end{equation}
It is the natural generalization of the usual Poisson kernel in the complex disc. It is proven in \cite{AS} that this kernel, which is a priori defined only for almost every $\xi\in N(\infty)$ admits a continuous extension to $N\times N\times N(\infty)$ which is characterized by the following properties ($z\in N$ and $\xi\in N(\infty)$):
$$k(z,z;\xi)=1,\,\,\,\,\,\,\Delta k(z,.\,;\xi)=0\,\,\,\,\,\,\textrm{and}\,\,\,\,\,\,\lim_{z_2\to\xi}k(z_1,z_2;\xi')=0\,\,\,\,\,\,\textrm{when}\,\,\,\,\,\,\xi'\neq\xi.$$

Note finally that the group $\Isom^+(N)$ of direct isometries acts diagonally on $N\times N\times\ N(\infty)$ and that the Poisson kernel satisfies the following property of equivariance which holds true for every $(z_1,z_2,\xi)\in N\times N\times N(\infty)$ and every $\gamma\in\Isom^+(N)$:
\begin{equation}
\label{equivariancekernel}
k(\gamma z_1,\gamma z_2;\gamma\xi)=k(z_1,z_2;\xi).
\end{equation}

In \cite{AS}, Anderson and Schoen study very precisely this kernel and prove the two following key theorems. First, the Poisson kernel is Hölder continuous.
\begin{theorem}[Hölder continuity of the kernel]
\label{poissonholder}
Fix two points $o,z\in N$. Then there exist positive constants $C>0$, $\alpha>0$ such that for any $\xi, \xi'\in N(\infty)$, 
$$|k(o,z;\xi)-k(o,z;\xi')|\leq C \angle_o(\xi,\xi')^{\alpha}.$$
\end{theorem}
Hence, since $k$ depends smoothly in the first two variables we have a Hölder dependance of $k$ with uniform constants in any $K_1\times K_2\times N(\infty)$ where $K_i$ are compact subsets of $N$. The second result gives an integral representation, also called \emph{Poisson representation}, of positive harmonic functions on $N$.

\begin{theorem}[Integral representation]
\label{integral}
Let $N$ be a complete connected and simply connected Riemannian manifold whose sectional curvature is pinched between two negative constants. Let $h$ be a positive harmonic function on $N$ and $o\in N$ be fixed. Then, there is a unique finite Borel measure $\eta_o$ on $N(\infty)$ such that for every $z\in N$, 
$$h(z)=\int_{N(\infty)}k(o,z;\xi)\,d\eta_o(\xi).$$
\end{theorem}

\paragraph{Remark.} This definition is independent of the choice of the base point $o$. Since the cocycle relation $k(o,o';\xi)k(o',z;\xi)=k(o,z;\xi)$ holds, the definition would also be valid with a base point $o'$ and a measure $\eta_{o'}=k(o,o';\xi)\eta_o$.

\subsection{Canonical lifts of harmonic measures}

We treat in this paragraph the first step in the proof of the main Theorem. Indeed, we show how to lift canonically any harmonic measure to the unit tangent bundle of a lamination whose leaves are negatively curved.

In what follows, $(M,\F)$ is a compact foliated space endowed with a leafwise metric such that the leaves are all negatively curved manifolds. By the compactness of $M$ the sectional curvatures of the leaves are uniformly pinched between two negative constants $-b^2\leq-a^2$. Also, the injectivity radii are uniformly bounded from below so we can choose a good foliated atlas $\A=(U_i,\phi_i)$ such that all plaques trivialize the universal cover.

\paragraph{Trivialization of the center unstable foliation.} Let $N$ be a complete, connected and simply connected Riemannian manifold whose sectional curvature is pinched between two negative constants $-b^2\leq-a^2$. The geodesic flow on $T^1N$ leaves invariant two continuous foliations. They are called the \emph{stable and unstable horospheric foliations} and are respectively denoted by $\widetilde{\W}^s$ and $\widetilde{\W}^u$. Their saturations by the flow are respectively called \emph{center stable} and \emph{center unstable} foliations and denoted by $\widetilde{\W}^{cs}$ and $\widetilde{\W}^{cu}$.

We have a map $\widetilde{\sigma}^+:T^1N\to N\times N(\infty)$ which associates to $v\in T^1N$ the couple $(c_v(0),c_v(-\infty))$, where $c_v$ is the geodesic determined by $v$. This map trivializes the center unstable foliation: a slice $N\times\{\xi\}$, $\xi\in N(\infty)$  is identifed the center unstable leaf associated to $\xi$ i.e. the set of vectors $v$ such that $c_v(-\infty)=\xi$.

Furthermore, this identification conjugates the actions of the group of direct isometries $\Isom^+(N)$ by differentials on $T^1N$ and by diagonal maps on $N\times N(\infty)$.

\paragraph{$\infty$-harmonic measures.} Fix a base point $o\in N$. The center unstable foliation $\widetilde{\W}^{cu}$ is identified with the trivial foliation $(N\times\{\xi\})_{\xi\in\widetilde{L}(\infty)}$. Thus it is possible to associate to each center unstable leaf $N\times\{\xi_0\}$ the $\xi_0$-\emph{harmonic measure} for $\widetilde{\W}^{cu}$ defined by:
$$m_{\xi_0}=k(o,z;\xi_0)\Leb_{N\times\{\xi_0\}}(z).$$

\begin{defi}
\label{inftyharmonic1}
A measure $m^+$ on $T^1N$ is said to be $\infty$-harmonic for $\widetilde{\W}^{cu}$ if it has a disintegration in the center unstable manifolds, identified with $(N\times\{\xi\})_{\xi\in N(\infty)}$, whose conditional measures are given by the $m_{\xi}$, $\xi\in N(\infty)$.

If $L$ is a quotient of $N$ by a subgroup of direct isometries, the quotient measure of an $\infty$-harmonic measure for $\widetilde{\W}^{cu}$ on $T^1L$ will also be called $\infty$-harmonic for the quotient center unstable foliation.
\end{defi}

\paragraph{Remark.} In the second part of the definition we used the fact that the measures $m_{\xi}$ satisfy the equivariance relation $\gamma\ast m_{\xi}=m_{\gamma\xi}$ for every direct isometry $\gamma$ (with the action on $N\times N(\infty)$ corresponding to the differential action on $T^1 N$). This comes from the fact that $\gamma$ acts as an isometry on center unstable leaves, and thus preserves the Lebesgue measure, as well as from the equivariance relation (\ref{equivariancekernel}).

\begin{lemma}
\label{projinfharmonic}
The projection on $N$ of any $\infty$-harmonic measure on $T^1N$ is a harmonic measure: it has a harmonic density with respect to Lebesgue.
\end{lemma}

\begin{demo}
Consider an $\infty$-harmonic measure $\widetilde{m}^+$ which is obtained by integration of measures $m_{\xi}$ against a finite measure $\eta_o$ on $N(\infty)$. All the center unstable leaves are isometric to $N$ and all $m_{\xi}$ have denstity with respect to Lebesgue. We deduce that the projection of $\widetilde{m}^+$ has a density with respect to Lebesgue.

Moreover, this density is obtained by integration of the $k(o,.;\xi)$ against $\eta_o$. Hence this is a harmonic function and the lemma is proven.
\end{demo}

\paragraph{Canonical lift.}  Remember that if $(M,\F)$ is a compact foliated space whose leaves are negatively curved the unit tangent bundle $T^1\F$ is foliated by the unit tangent bundles of the leaves of $\F$. This lamination is denoted by $\hcF$ and is itself subfoliated by the center unstable foliation $\W^{cu}$ of the foliated geodesic flow $G_t$.

\begin{defi}
Let $(M,\F)$ be a compact foliated space with a leafwise metric. Assume that all the leaves of $\F$ are negatively curved. A probability measure $m^+$ on $T^1\F$ will be called $\infty$-harmonic for $\W^{cu}$ if its conditional measures on the plaques of $\hcF$ are proportional to the $\infty$-harmonic for the center unstable foliations of the corresponding leaves (see Definition \ref{inftyharmonic1} for the definition of $\infty$-harmonic measures on a leaf $L$).
\end{defi}

Before we state the next proposition let us introduce three notations:
\begin{itemize}
\item $pr:T^1\F\to M$ stands for the canonical projection along the tangent spheres;
\item $\HH ar(\F)$ stands for the set of harmonic measure for $\F$;
\item $\HH ar_{\infty}(\W^{cu})$ stands for the set of $\infty$-harmonic measures for $\W^{cu}$.
\end{itemize}

In the next proposition, we present an ``unrolling argument'' allowing us to lift canonically the harmonic measures to $T^1\F$. We will obtain a section $\HH ar(\F)\to\HH ar_{\infty}(\W^{cu})$.

\begin{proposition}
\label{liftharmonicmeasures}
Let $(M,\F)$ be a compact foliated space with a leafwise metric such that all the leaves are negatively curved. Then:
\begin{enumerate}
\item if $m^+$ is an $\infty$-harmonic measure for $\W^{cu}$, the projection $pr_{\ast} m^+$ is a harmonic measure for $\F$;
\item there is an injective map $\sigma:\HH ar(\F)\to\HH ar_{\infty}(\W^{cu})$ which satisfies for any harmonic measure $m$ $pr_{\ast}[\sigma(m)]=m$. The measure $\sigma(m)$ is called the canonical lift of $m$.
\end{enumerate}
In other terms, the map $pr_{\ast}:\HH ar_{\infty}(\W^{cu})\to\HH ar(\F)$ which associates to any $\infty$-harmonic measure $m^+$ its projection on $M$, $pr_{\ast} m^+$ is well defined, surjective and has a section given by the canonical lift.
\end{proposition}

\begin{demo}
Choose a good foliated atlas $\A=(U_i,\phi_i)_{i\in I}$ for $\F$ such that when we pull it back by the fibration we obtain a good foliated atlas for $\hcF$ whose charts $\widehat{U}_i$ are trivially foliated by the unit tangent spheres and trivialize $\W^{cu}$. We can also assume that the plaques of $U_i$ trivialize the Riemannian cover of every leaf.

The fact that the function $pr_{\ast}$ is well defined, or if one prefers that the projection on $M$ of an $\infty$-harmonic measure for $\W^{cu}$ is harmonic, follows from Lemma \ref{projinfharmonic} applied to the plaques of $\F$.

Now let us show how to lift a harmonic measure to the unit tangent bundle. Because all harmonic measures are barycenters of ergodic components (cf \cite{Gar}), we only have to show how to lift ergodic harmonic measures. So let $m$ be an ergodic harmonic measure. By ergodicity, there exists some $i_0\in I$ such that the saturation of $T_{i_0}$ by $\F$ has full measure. In other terms, almost every leaf intersects the transversal $T_{i_0}$.

Thus, by Lemma \ref{harmonicextension}, there is a saturated Borel set $\X\dans M$ which is full for $m$ with the following properties: all leaves passing through a point of $\X$ intersect $T_{i_0}$ and for any $y\in\X$ the harmonic density defined on a corresponding plaque can be extended to the whole leaf.

Let $x_0\in T_{i_0}\cap\X$. Formula (\ref{Eq:harmonicextension1}) defines a harmonic function $H_{x_0}:L_{x_0}\to(0,\infty)$. We can lift it to the Riemannian universal cover $\widetilde{L}_{x_0}$ (fix a preimage $o$ of $x_0$) as a harmonic function $\widetilde{H}_{x_0}$.

Since the plaques $P_i$ trivialize the universal cover, any $h_i$ can be lifted as a harmonic function $\widetilde{h_i}$ of any section $\widetilde{P_i}\dans\widetilde{L}_{x_0}$. By Lemma \ref{harmonicextension}, if $c$ is a path joining $x_0$ and the plaque $P_i$ and if $\widetilde{P_i}(c)$ is the corresponding lift, then the formula $\widetilde{H}_{x_0}=\widetilde{h_i}d[\tau_c^{-1}\ast\nu_i]/d\nu_{i_0}(x_0)$ holds in restriction to $\widetilde{P_i}(c)$.

Theorem \ref{integral} gives us a measure $\eta_o$ on $\widetilde{L}_{x_0}(\infty)$ such that:
$$\widetilde{H}_{x_0}(z)=\int_{L_{x_0}(\infty)}k(o,z;\xi) d\eta_{o}(\xi).$$
Then, it is possible to lift the measure $\widetilde{H}_{x_0}\,\Leb$ to $T^1\widetilde{L}_{x_0}$ by considering the measure $\widetilde{m}^+_{x_0}$ obtained by integration against $\eta_o$ of the measures $m_{\xi}$ defined on the unstable leaves identified with $\widetilde{L}_{x_0}\times\{\xi\}$. Recall that we have an equivariance relation $\gamma\ast m_{\xi}=m_{\gamma\xi}$ which holds for every $\gamma\in\pi_1(L_{x_0},x_0)$ and that $\widetilde{H}_{x_0}$ is invariant by the action of $\pi_1(L_{x_0},x_0)$. Hence the measure $\widetilde{m}^+_{x_0}$ is invariant under the action of $\pi_1(L,x_0)$ and by projection onto $T^1L_{x_0}$, we obtain a measure defined on the whole unit tangent bundle of $L_{x_0}$. We denote it by $m^+_{x_0}$ and if we push it by $pr$ it projects down onto $H_{x_0}\,\Leb$.

Now consider $y\in L_{x_0}\cap\TT$. There exists an index $i\in I$ such that $y\in T_i$. One can form a measure $m^+_{i,y}$ on the plaque $T^1P_i(y)$ by dividing the restriction of $m^+_{x_0}$ to $T^1P_i(y)$ by the derivative $d[\tau_c^{-1} \nu_i]/d\nu_{i_0}(x_0)$, where $c$ is any path joining $x_0$ to $y$.

Another application of Lemma \ref{harmonicextension} shows that this measure does not depend on the base point $x_0$ Indeed if $x_1\in T_{i_0}\cap L_{x_0}$, we have:
$$\left(\frac{d[\tau_c^{-1}\ast\nu_i]}{d\nu_{i_0}}(x_0)\right)^{-1}(m^+_{x_0})_{|T^1P_i(y)}=\left(\frac{d[\tau_{c'}^{-1}\ast\nu_i]}{d\nu_{i_0}}(x_1)\right)^{-1}(m^+_{x_1})_{|T^1P_i(y)},$$
where $c$ and $c'$ are any paths joining respectively $x_0$ and $x_1$ to $y$. Moreover, when we push the measure $m_{i,y}$ on the plaque $P_i(y)$, it projects down onto $h_i(.,y)\,\Leb_{P_i(y)}$ .

Assume that $m(U_i\cap U_j)>0$ and that $y\in\dom(\tau_{ij})$, then we have for any $v\in T^1P_i(y)\cap T^1P_j(\tau_{ij}(y))$,
\begin{equation}
\label{Eq:gluemeasurestogether}
\frac{dm^+_{i,y}}{dm^+_{j,\tau_{ij}(y)}}(v)=\frac{d[\tau_{ij}^{-1}\ast\nu_j]}{d\nu_i}(y),
\end{equation}
in such a way that we will be able to glue these measures together.

Finally, almost every leaf meets $T_{i_0}$. So the construction above gives us for $\widehat{m}$-almost every $y\in \TT$ a measure $m_{i,y}^+$ on $T^1P_i(y)$ (if $y\in T_i$) which depends only on $i$ and on $y$. Note that this family of measures is $\nu_i$-measurable. Indeed it comes directly from the fact the family of local densities $h_i(.,y)$ varies measurably with $y\in T_i$: hence the lift $(d[\tau_c^{-1}\ast\nu_i]/d\nu_{i_0}(x_0))^{-1} \widetilde{H}_{x_0}$ varies measurably with $y$. The correspondence $\widetilde{H}_{x_0}\mapsto \eta_o$ being a homeomorphism between harmonic functions which are $1$ at $o$ and probability measures on the sphere at infinity, we have the measurability with $y$ of the family $(d[\tau_c^{-1}\ast\nu_i]/d\nu_{i_0}(x_0))^{-1}\widetilde{m}^+_{x_0}$, which is exactly the lift of the family $m^+_{i,y}$.

Hence, for any $i$ such that $m(U_i)>0$ we obtain a measure $m_i$ on $\widehat{U}_i$ by integrating the measures $m_{i,y}^+$ against $\nu_i$. Because of Relation (\ref{Eq:gluemeasurestogether}), we see that these measures can be glued together. Hence, we have a well defined measure $\sigma(m)=m^+$ on $T^1\F$.

By definition, the measure induced by $m^+$ on $\TT$ is exactly given by $\widehat{m}$. Furthermore, since each of the conditional measures of $m^+$ in the plaques are $\infty$-harmonic and project down to that of $m$, we conclude that $m^+$ is $\infty$-harmonic for $\W^{cu}$ and projects down to $m$.
\end{demo}

\paragraph{Remark.} In order to prove that the map $pr_{\ast}:\HH ar_{\infty}(\W^{cu})\to\HH ar(\F)$ is a bijection, it remains to prove that it is injective. A proof of this fact might be adapted from the proof of \cite{BMar} in the case of lamination by hyperbolic Riemann surfaces. It would rely on an analytic argument based on the uniqueness of the integral representation (see Theorem \ref{integral}). We propose another proof more in the spirit of this work. We will prove that two different $\infty$-harmonic measures for $\W^{cu}$ induce different measures on the transversal of $\hcF$: hence they project onto different harmonic measures. This is based on a variation of the Hopf argument, but we first need the definition of $H$-Gibbs measures in order to explain this fact.

\section{$H$-Gibbs measures for the foliated geodesic flow}
Now that we know how to lift canonically harmonic measures to the unit tangent bundle, the next step is to reparametrize these lifts in the center unstable direction so as to obtain invariant measures with the Gibbs property.

\subsection{Harmonic class on horospheres}
We assume here that $N$ is a complete, connected and simply connected Riemannian manifold whose sectional curvature is pinched between two negative constants.

\paragraph{Measures on strong and weak unstable leaves.} We attribute the next result to Ledrappier \cite{L2}. It gives a family of measures in the unstable leaves satisfying a relation of quasi-invariance by the geodesic flow of Gibbs type associated to the following potential:
$$H(v)=\left.\frac{d}{dt}\right|_{t=0}\log\,k(c_v(0),c_v(t);c_v(-\infty)),$$
where $c_v$ is the geodesic directed by $v$. Note that we don't require $N$ to be a Riemannian cover of a compact manifold, but only the pinching of sectional curvatures. By projection of the harmonic class on the horospheres and multiplication by a suitable function, Ledrappier was able to prove the following result:

\begin{theorem}
\label{harmonicunstable}
Let $N$ be a complete connected and simply connected Riemannian manifold whose sectional curvature is pinched between two negative constants. There exist two families of measures denoted by $(\widetilde{\lambda}^{u}_{H,v})_{v\in T^1N}$ and $(\widetilde{\lambda}^{cu}_{H,\xi})_{\xi\in N(\infty)}$ respectively defined on the strong and weak unstable leaves of the geodesic flow $G_t$ such that:
\begin{enumerate}
\item for every $v\in T^1N$, $\widetilde{\lambda}^u_{H,v}$ is equivalent to the projection of the harmonic class on $\widetilde{W}^u(v)$ along the geodesics;
\item $\widetilde{\lambda}^{u}_{H,v}$ is independent of the choice of $v$ in a given unstable leaf;
\item we have the following equivariance relations for any $\gamma\in\Isom^+(N)$, $v\in T^1N$ and $\xi\in N(\infty)$: 
$$\gamma\ast\widetilde{\lambda}^{u}_{H,v}=\widetilde{\lambda}^{u}_{H,D\gamma v}\,\text{   and   }\,\gamma\ast\widetilde{\lambda}^{cu}_{H,\xi}=\widetilde{\lambda}^{cu}_{H,\gamma\xi};$$
\item we have the following quasi-invariance property by the geodesic flow:
\begin{equation}
\label{Eq:harmoniccocycle1}
\frac{d\left[G_T\ast\widetilde{\lambda}^u_{H,G_{-T}(v)}\right]}{d\widetilde{\lambda}^u_{H,v}}(w)=k(c_w(-T),c_w(0);c_w(-\infty))
\end{equation}
where $T\in\R$, $v\in T^1N$, $w\in W^u(v)$ and $c_w$ denotes the geodesic directed by $w$;
\item for any $v\in T^1N$ with $c_v(-\infty)=\xi$, the measure $\widetilde{\lambda}^{cu}_{H,\xi}$ is obtained by integration of the $\widetilde{\lambda}^u_{H,G_t(v)}$ against the $dt$ element of the flow: in particular, this measure satisfies Relation (\ref{Eq:harmoniccocycle1}).
\end{enumerate}
\end{theorem}

\paragraph{Absolute continuity.} We will need the following absolute continuity property of the stable foliation.

\begin{theorem}
\label{absolutecontinuityharmonic}
Let $N$ be a complete connected and simply connected Riemannian manifold whose sectional curvature is pinched between two negative constants. Then the holonomy maps along the stable foliation preserve the class of $\lambda^{cu}_{H,\xi}$. More precisely, let for $v,w\in T^1N$ in the same stable leaf:
\begin{equation}
\label{Eq:jacobian}
k^s_H(v,w)=\exp\left[\int_0^{\infty}(H\circ G_t(w)-H\circ G_t(v))dt\right]=\lim_{T\to\infty}\frac{k(c_w(0),c_w(T);c_w(-\infty))}{k(c_v(0),c_v(T);c_v(-\infty))}.
\end{equation}
Then:
\begin{enumerate}
\item $k^s_H(v,w)$ is well defined for all $v,w$ lying in the same stable manifold;
\item for all $v_0$ and $w_0$ lying on the same stable manifold, one has, for any $w\in W^{cu}(w_0)$,
$$\frac{d\left[\h s{v_0}{w_0}\ast\widetilde{\lambda}^{cu}_{H,v_0}\right]}{d\lambda^{cu}_{H,w_0}}(w)=k^s_H(w,\h s{w_0}{v_0}(w));$$
\item there exist uniform Hölder constants $C_0$ and $\alpha_0$ such that when $v$ and $w$ lie in the same stable leaf with $\dist_s(v,w)\leq 1$, $\log k^s_H(v,w)\leq C_0\dist_s(v,w)^{\alpha_0}$;
\end{enumerate}
\end{theorem}

\begin{demo}
First, the integral given in Formula (\ref{Eq:jacobian}) exists by a usual distortion lemma, because the potential $H$ is Hölder continuous with uniform constants (which depends only on the bounds of the curvature): see Theorem \ref{poissonholder}.

Second, the absolute continuity has been proven in \cite{L2}: it only relies on the pinching of the curvature. Indeed, since the Anosov constants and the Hölder constants for the kernel depend only on the pinching of the curvature (see \cite{AS}), the distortion controls that are usually needed to deduce absolute continuity (see for example Lemmas 3.3 and 3.6 of Ma\~né's book \cite{Ma}) are proven in a similar way. In particular, the formula for the Jacobian $k^s_H$ is given by Formula (\ref{Eq:jacobian}) and we have the uniform log-Hölder continuity of the Jacobian $k^s_H$ (see the end of the proof of Theorem 3.1 in \cite{Ma}).
\end{demo}

\subsection{$H$-Gibbs measures for the foliated geodesic flow}

Until the end of this article, $(M,\F)$ will denote a compact foliated space endowed with a foliated metric such that all of its leaves are negatively curved.

\paragraph{$H$-Gibbs measures.} We know that there is a family of measures denoted by $(\lambda^u_{H,v})_{v\in T^1\F}$ on the unstable leaves of $G_t$ which satisfy the cocycle Relations (\ref{Eq:harmoniccocycle1}). There is also a family of measures defined on the center unstable leaves, denoted by $(\lambda^{cu}_{H,v})_{v\in T^1\F}$, which is obtained by integration of the $\lambda^u_{H,G_t(v)}$ against the $dt$ element of the flow.

\begin{defi}
\label{defHgibbs}
An $H$-Gibbs measure for the foliated geodesic flow $G_t$ is a probability measure $\mu$ on $T^1\F$, which satisfies the following properties:
\begin{itemize}
\item $\mu$ is invariant by $G_t$;
\item $\mu$ has an absolutely continuous disintegration with respect to $(\W^u,\lambda^u_{H,v})$.
\end{itemize}
\end{defi}

The letter $H$ stands for harmonic. Also, we would like to recall that these measures are associated to the following potential:
$$H(v)=\left.\frac{d}{dt}\right|_{t=0}\log\,k_{L(v)}(c_v(0),c_v(t);c_v(-\infty)),$$
where $L(v)$ is the leaf $v$ is tangent to, $k_{L(v)}$ is the Poisson kernel on $\widetilde{L(v)}$ and $c_v$ is the geodesic directed by $v$.

\begin{theorem}
\label{existenceofhgibbs}
Let $(M,\F)$ be a compact foliated space with a leafwise metric. Assume that the leaves have negative sectional curvature.
\begin{enumerate}
\item For any $v\in T^1\F$ and any small disc $D\dans\W^u_{loc}(v)$, the accumulation points of
$$\frac{1}{T}\int_0^T\frac{G_t\ast(\lambda^u_{H,v})_{|D}}{\lambda^u_{H,v}(D)} dt$$
are $H$-Gibbs measures for $G_t$.
\item Ergodic components of $H$-Gibbs measures for $G_t$ are $H$-Gibbs measures.
\item If $\mu$ is any $H$-Gibbs measure for $G_t$, then the local densities $\psi^u_{H,v}$ on the local unstable manifolds are uniformly log-bounded and satisfy, for $w\in W^u_{loc}(v)$:
\begin{equation}
\label{Eq:characgibbsmeasures1}
\frac{\psi^u_{H,v_0}(v)}{\psi^u_{H,v_0}(w)}=k_{L(v)}(c_w(0),c_v(0);c_v(-\infty)).
\end{equation}
\end{enumerate}
\end{theorem}

\begin{demo}
The proof of this theorem is a simple adaptation of the section 11.2.2. of \cite{BDV}, where the authors study Gibbs u-states of partially hyperbolic systems. Here, we know that $\W^u$ is a lamination of a \emph{compact} space $M$ which is \emph{uniformly} expanded by the flow $G_t$. By using Formula (\ref{Eq:harmoniccocycle1}) and the uniform Hölder continuity of the log of the Poisson kernel on the leaves (see Theorem \ref{poissonholder}), we see that these results can be generalized without difficulty (the precise verifications can be found in the author's thesis \cite{Al3}).
\end{demo}

\paragraph{Remark.} In the definition of $H$-Gibbs measure, we have been led to make a choice. Indeed, we could have also defined the notion of Gibbs measure by prescribing the class of conditional measures in the stable manifolds. The important fact is that in this foliated case we can't prescribe the classes of conditional measures in the stable and unstable manifolds simultaneously. An illustration of this fact is a result proved in \cite{M}: for a lamination by hyperbolic Riemann surfaces, the existence of measures invariant by the actions of the foliated stable and unstable horocyclic flows implies that of a transverse measure invariant by holonomy, which is very rare.

\begin{proposition}
\label{transverseGibbs}
Let $(M,\F)$ be a compact foliated space with a leafwise metric. Assume that the leaves have negative sectional curvature. Let $\B=(V_i,g_i)_{i\in I}$ be a good foliated atlas for $\bcW^u$ in $T^1\F$. Denote by $(S_i)_{i\in I}$ the associated complete system of transversals. Let $\mu$ be an $H$-Gibbs measure. Then, there is a family of measures $(\nu_i)_{i\in I}$ on these transversals such that:
\begin{enumerate}
\item $\nu_i$ is equivalent to the projection on $T_i$ of the restriction $\mu_{|V_i}$ with a Radon-Nikodym derivative which is log-bounded independently of $i$;
\item it satisfies the following quasi-invariance relation:
\begin{equation}
\label{Eq:harmoniccocycle2}
\frac{d[\h u{S_i}{S_j}\ast\nu_i]}{d\nu_j}(w)=k_{L(v)}(c_w(0),c_v(0);c_v(-\infty))
\end{equation}
for $v$ in the domain of $\h u{S_j}{S_i}$ and $w=\h u{S_j}{S_i}$.
\end{enumerate}
\end{proposition}

\begin{demo}
By Theorem \ref{existenceofhgibbs} the conditional measure of an $H$-Gibbs measure $\mu$ for $G_t$ in the unstable plaques $W^u_{loc}(v)$ have the form $\psi^u_{H,v}\,\bar{\lambda}^u_{F,v}$, $v_0\in S_i$.

We are free to choose $\psi^u_{H,v}$ equal to $1$ in each $S_i$: by uniqueness of the disintegration, it determines the transverse measure $\nu_i$. Because of the Hölder continuity of the Poisson kernel, such a measure obviously has a log-bounded derivative with respect to the projection on $S_i$ of the restriction $\mu_{|V_i}$

Since the densities satisfy relations (\ref{Eq:characgibbsmeasures1}), we obtain the desired cocycle relation by evaluating $\mu$ on an intersection $V_i\cap V_j$.
\end{demo}

\subsection{Induced measures on a complete transversal} It is true that two different $H$-Gibbs measures induce different measures on a complete transversal to $\W^u$ (see Proposition \ref{transverseGibbs}). It is more interesting to note, and this is the goal of the following proposition, that an $H$-Gibbs measure is in fact determined by the measure it induces on a complete transversal to $\hcF$.

In order to prove this fact, we propose an argument ``à la Hopf'', where the principal ingredient is the absolute continuity property \ref{absolutecontinuityharmonic} and that we will generalize in a further study of Gibbs measures for the foliated geodesic flow (see \cite{Al1,Al2}). We have chosen a good foliated atlas $\A$ of $\F$. We denote by $\TT$ a corresponding complete transversal. Recall that by Lemma \ref{sameholonomy} by pulling back $\A$ by the fibration we obtain a foliated atlas $\widehat{A}$ for $\hcF$ whose complete transversal can be identified with $\TT$: the pseudogroups of $\hcF$ and $\F$ are identified.

\begin{proposition}
\label{mesuresinduiteshgibbs}
Let $(M,\F)$ be a compact foliated space with a leafwise metric. Assume that the leaves have negative sectional curvature. Let $\TT$ be a complete transversal to $\F$, identified with a complete transversal of $\hcF$. Then, we have that:
\begin{enumerate}
\item any $H$-Gibbs measure for $G_t$ induces a measure on $\TT$ which is quasi-invariant by the holonomy pseudogroup of $\hcF$;
\item two different $H$-Gibbs measures for $G_t$ induce different measures on $\TT$.
\end{enumerate}
\end{proposition}

\begin{demo}
Let $\mu$ be an $H$-Gibbs measure for $G_t$. Let $\widehat{U}_i$ be the charts of $\widehat{\A}$ and let $(T_i)_{i\in I}$ be the complete system of transversal: we have $\TT=\bigsqcup_{i\in I}T_i$. Consequently, one obtains a complete system of transversals for $\W^u$ by considering the $T^{cs}_i=\bigcup_{x\in T_i}W^{cs}_{loc}(x)$. By Proposition \ref{transverseGibbs}, if $\nu^{cs}_i$ denotes the projection of $\mu_{|U_i}$ on the transversals $T_i^{cs}$, the family $(\nu^{cs}_i)_{i\in I}$ is quasi-invariant by the holonomy maps along the unstable leaves.

Denote by $\nu_i$ the projection of $\mu_{|\widehat{U}_i}$ on $T_i$ along the leaves of $\F$. It coincides with the projection of $\nu_i^{cs}$ along the center stable leaves. Now assume that two charts $\widehat{U}_i$ and $\widehat{U}_j$ have an intersection of positive measure for $\mu$, or, it amounts to the same thing, that the domain of $\tau_{ij}$ is positive for $\nu_i$. If $A$ lies in the domain of $\tau_{ij}$ and satisfies $\nu_i(A)>0$, then we have $\nu^{cs}(\overline{A})>0$ where $\overline{A}$ is the saturated, inside $T^{cs}_i$, of $A$ in the center stable direction. As we have said before, $\nu^{cs}_j(\h u{T^{cs}_i}{T^{cs}_j}(\overline{A}))>0$. By projection along the center stable direction, we have $\nu_j(\tau_{ij}(A))>0$. This proves that $\nu_i$ and $\tau_{ji}\ast\nu_j$ are in the same measure class. By induction, we prove that the holonomy pseudogroup (generated by the $\tau_{ij}$) preserves this measure class.\\

Now we have to prove that the classes of measure induced on $\TT$ by two different $H$-Gibbs measures for $G_t$ are different. Since any ergodic component of an $H$-Gibbs measure for $G_t$ is still an $H$-Gibbs measure, this is the same as proving that two different (hence singular) ergodic $H$-Gibbs measures for $G_t$ induce singular measures on $\TT$.

This is where we have to use the absolute continuity property \ref{absolutecontinuityharmonic}. Consider two ergodic $H$-Gibbs measures $\mu_1$ and $\mu_2$ and assume that the induced measures on $\TT$, denoted by $\widehat{\mu}_1$ and $\widehat{\mu}_2$, aren't singular. Let $\X$ be the set of vectors $v\in T^1\F$ such that for any continuous function $f:T^1\F\to\R$ we have:
$$\lim_{T\to\infty}\frac{1}{T}\int_0^Tf\circ G_t(v)dt=\int_{T^1\F} fd\mu_1.$$
By Birkhoff's ergodic theorem, we know that $\mu_1(\X)=1$. Hence, if we denote by $\widehat{X}_i$ the projections of $\X\cap \widehat{U}_i$ on transversal $T_i$ and $\widehat{X}=\bigsqcup_i\widehat{X}_i$ we have by definition, $\widehat{\mu}_1(^c\widehat{X})=0$. Since we assumed that $\widehat{\mu}_1$ and $\widehat{\mu}_2$ aren't singular we have $\widehat{\mu}_2(\widehat{X})>0$.

Hence there exists $i\in I$ that satisfies $\widehat{\mu}_2(\widehat{X}_i)>0$ where $\widehat{X}_i=\widehat{X}\cap T_i$. Thus there is a set $\Y\dans\widehat{U}_i$ which is positive for $\mu_2$ and projects down onto $\widehat{X}_i$ along the plaques. By ergodicity of $\mu_2$, we can assume that for any $w\in\Y$ and for any continuous function $f:T^1\F\to\R$ we have:
$$\lim_{T\to\infty}\frac{1}{T}\int_0^Tf\circ G_t(w)dt=\int_{T^1\F} fd\mu_2.$$

Now consider a point of $\widehat{X}_i$, that we will identify with a plaque of $\hcF$ denoted by $\widehat{P}_i$. Since $\mu_1$ and $\mu_2$ have an absolutely continuous disintegration with respect to $(\W^{cu},\lambda^{cu}_{H,v})_{v\in T^1\F}$, we can find two vectors $v\in\X\cap \widehat{P}_i$ and $w\in\Y\cap \widehat{P}_i$ such that $\lambda^{cu}_{H,v}$-almost every vector of $W^{cu}_{loc}(v)$ belongs to $\X$ and that $W^{cu}_{loc}(w)$ meets $\Y$ on a Borel set of positive $\lambda^{cu}_{H,w}$-measure.

Let us project the intersection $\Y\cap W^{cu}_{loc}(w)$ on $W^{cu}_{loc}(v)$ along the stable foliation. By the property of absolute continuity \ref{absolutecontinuityharmonic}, this projection has positive $\lambda^{cu}_{H,v}$ measure in $W^{cu}_{loc}(v)$ and thus intersects $\X$. To summarize, we proved the existence of $v'\in\X$ and $w'\in\Y$ which are on the same local stable manifold. Hence, for any continuous function $f:T^1\F\to\R$:
$$\lim_{T\to\infty}\frac{1}{T}\int_0^Tf\circ G_t(v')dt=\lim_{T\to\infty}\frac{1}{T}\int_0^Tf\circ G_t(w')dt.$$
But the first limit equals $\int_{T^1\F} fd\mu_1$ and the second one equals $\int_{T^1\F} fd\mu_2$: these two integrals must be the same, and this for any continuous function $f:M\to\R$. The equality of the two measures follows.
\end{demo}

\paragraph{Conditional measures on stable manifolds.} It is also possible to disintegrate an $H$-Gibbs measure $\mu$ for $G_t$ in the stable foliation.
\begin{enumerate}
\item By invariance by the foliated geodesic flow and by Relation (\ref{Eq:harmoniccocycle1}) the conditional measures $\lambda^s_{H,v}$ on \emph{local} stable manifold of the $H$-Gibbs measure $\mu$ satisfy the following quasi-invariance Relation:
\begin{equation}
\label{Eq:harmoniccocycle3}
\frac{d\left[G_T\ast\lambda^s_{H,G_{-T}(v)}\right]}{d\lambda^s_{H,v}}(w)=k_{L(v)}(c_w(0),c_w(-T);c_w(-\infty))
\end{equation}
Hence associated to each $H$-Gibbs measure we have a well defined class of conditional measures in \emph{global} stable manifolds which is quasi-invariant by the foliated geodesic flow.
\item This class is essentially unique in the sense that if there is class of measure on local stable manifold satisfying the cocycle Relation (\ref{Eq:harmoniccocycle3}) it is possible to combine it with $\lambda^{cu}_{F,v}$ in order to get an $H$-Gibbs measure which induces the same transverse measure on a complete transversal as $\mu$: these measures have to be equal by Proposition \ref{mesuresinduiteshgibbs}.
\item By Proposition \ref{transverseGibbs} if we integrate the family of $\lambda^s_{H,v}$ against the $dt$-element we obtain a family of measures on center stable manifolds which is invariant by unstable holonomies with cocycle Relations (\ref{Eq:harmoniccocycle2}).
\item It is an interesting problem to identify this class of measure. We conjecture that if it is equivalent to the harmonic class on stable manifold then there exists a tranvserse invariant measure. This conjecture is inspired by the main result of \cite{Al1} where it is shown that if an invariant measure has Lebesgue disintegration both in the stable and unstable leaves, then there exists a transverse invariant measure. Also this conjecture is true in the context of foliated $\C\PP^1$-bundles \cite{Al2}.
\end{enumerate}

\subsection{Bijective correspondence between harmonic measures and $H$-Gibbs measures}

\paragraph{Reparametrization of the canonical lift.} 

In Proposition \ref{liftharmonicmeasures}, we defined the canonical lift $m^+$ of any harmonic measure $m$. We will associate to this measure an $H$-Gibbs measure by reparametrisation of $m^+$ along the center unstable manifolds. We will exchange the Lebesgue measure on center unstable leaves by the measure $\lambda^{cu}_{H,v}$.

\begin{lemma}
\label{reparamxi}
Let $L$ be a leaf of $\F$ and $\xi_0\in\widetilde{L}(\infty)$. The following measure defined on the center unstable leaf identified to $\widetilde{L}\times\{\xi_0\}$ is invariant by the geodesic flow $G_t$, and has an absolutely continuous disintegration with respect to $(\widetilde{\W}^{u},\widetilde{\lambda}^{u}_{H,\xi})$:
$$\mu_{\xi_0}=k(o,z;\xi_0)\widetilde{\lambda}^{cu}_{H,\xi_0}(z).$$
\end{lemma}

\begin{demo}
The flow invariance is an immediate consequence of Formula (\ref{Eq:harmoniccocycle1}). Since $\widetilde{\lambda}^{cu}_{H,\xi}$ is obtained by integration of the $\widetilde{\lambda}^u_{H,G_t(v)}$ against the $dt$ element of the flow, this measure has an absolutely continuous disintegration with respect to $(\widetilde{\W}^{cu},\widetilde{\lambda}^{cu}_{H,\xi})$.
\end{demo}

\paragraph{Remark.} We also have the following equivariance relation which holds for every $\xi_0\in N(\infty)$ and $\gamma\in\pi_1(L)$ and follows from Relation \ref{equivariancekernel}:
$$\gamma\ast\mu_{\xi_0}=\mu_{\gamma\xi_0},$$
so that this family of measures can be projected down on $T^1L$.

\begin{proposition}
\label{reparamharmhgibbs}
Let $(M,\F)$ be a compact foliated space with a leafwise metric. Assume that the leaves of $\F$ are negatively curved. Let $\TT$ be a complete transversal. Then:
\begin{itemize}
\item for any $\infty$-harmonic measure for $\W^{cu}$ $m^+$, there exists a unique $H$-Gibbs measure $\mu$ which induces on $\TT$ the same measure as $m^+$;

\item reciprocally, for any $H$-Gibbs measure $\mu$, there exists a unique $\infty$-harmonic measure for $\W^{cu}$ $m^+$ which induces on $\TT$ the same measure as $\mu$.
\end{itemize}
\end{proposition}

\begin{demo}
Choose a foliated atlas $\widehat{\A}$ for $\hcF$ which trivializes the center unstable foliation. The center unstable plaques of this atlas are denoted by $Q_i$ and the transversals are denoted by $S_i$. The disjoint union of the $S_i$ gives a complete transversal to $\W^{cu}$ denoted by $\cS$. It can be seen as a trivial foliation in spheres $\bigsqcup_{x\in\TT} T^1_xL$, where $\TT$ is a complete transversal to $\F$.

The measure $m^+$ is by definition obtained in a chart $U_i$ by integration against a measure $\nu'_i$ defined on $S_i$ of measures which are given by the projection on the plaque $Q_i$ of the corresponding $m_{\xi}$.

Instead of integrating the measures $m_{\xi}$, let us integrate the projection on $Q_i$ of $\mu_{\xi}$ defined in Lemma \ref{reparamxi} (see the remark above for the equivariance of the family $(\mu_{\xi})_{\xi\in\widetilde{L}(\infty)}$) against $\nu_i'$ and denote by $\mu_i$ the measures we obtain on the charts $\widehat{U}_i$.

On the one hand, these measures $\mu_i$ are invariant by the foliated geodesic flow and have absolutely continuous disintegration with respect to $(\W^{u},\lambda^u_{H,v})$ because so does each $\mu_{\xi}$ (see Lemma \ref{reparamxi}).

On the other hand, for $v\in T^1\F$, the measure $\lambda^{cu}_{H,v}$ is, like the Lebesgue measure, well defined on all the center unstable leaf of $v$. Hence, like the measures $m^+_i$, the measures $\mu_i$ can be glued together. The resulting measure, denoted by $\mu$, is invariant by the flow $G_t$ and has an absolutely continuous disintegration with respect to $(\W^u,\lambda^u_{H,v})$. This is an $H$-Gibbs measure, and by construction, it induces the same transverse measure as $m^+$ on a complete transversal $\cS$. This implies two things.

First, since $\cS$ is a trivial foliation in spheres parametrized by $\TT$, $m^+$ and $\mu$ induce the same measure on $\TT$.

Second, since an $\infty$-harmonic measure is characterized by the measure it induces on $\cS$, two different $\infty$-harmonic measures yield two different $H$-Gibbs measures: we have the uniqueness.

In order to have the converse result, read the latter argument backwards.
\end{demo}

Thus, since we know that two different $H$-Gibbs measures induce on $\TT$ two different measures (see Proposition \ref{mesuresinduiteshgibbs}), it comes that this property also holds for two different $\infty$-harmonic measures for $\W^{cu}$. Then, we have the following proposition:
\begin{proposition}
\label{projinjective}
The map $pr_{\ast}:\HH ar_{\infty}(\W^{cu})\to\HH ar(\F)$ is injective.
\end{proposition}

\paragraph{End of the proof of the main Theorem.} We are now ready to prove the main Theorem.

On the one hand, Propositions \ref{liftharmonicmeasures} and \ref{projinjective} show that for any harmonic measure for $\F$ there is a unique $\infty$-harmonic measure for $\W^{cu}$ inducing the same measure on $\TT$, and reciprocally.

On the other hand, Proposition \ref{reparamharmhgibbs} shows that for any $\infty$-harmonic measure for $\W^{cu}$ there is a unique $H$-Gibbs measure for $G_t$ inducing the same measure on $\TT$ and reciprocally.

Hence, the proof of the theorem follows.\quad \hfill $\square$

\vspace{10pt}
\paragraph{Acknowledgments} This article corresponds to the fourth chapter of my PhD thesis \cite{Al3}. I am deeply grateful to my advisor Christian Bonatti for its constant patience and encouragements. I also wish to thank the anonymous referee whose comments helped to improve this paper. During the preparation of my master's thesis, Marco Brunella very kindly suggested me to write an unrolling argument in order to shorten the proof of the theorem of Mart\'inez and Bakhtin-Mart\'inez. I wish to dedicate this work to his memory.

\vspace{10pt}

\noindent \textbf{Sébastien Alvarez (sebastien.alvarez@u-bourgogne.fr)}\\
\noindent  Institut de Math\'ematiques de Bourgogne, CNRS-UMR 5584\\
\noindent  Université de Bourgogne, 21078 Dijon Cedex, France\\

\end{document}